\documentclass[centertags, reqno]{amsart}                
\usepackage{graphicx}
\usepackage{url}    
\usepackage{dsfont}

\newtheorem{thm}{Theorem}
\newtheorem*{thm*}{Theorem}

\theoremstyle{definition}

\theoremstyle{remark}

\newcommand{\mr}{{\mathbb R}}

\newcommand{\mn}{{\mathbb N}}

\newcommand{\mc}{{\mathbb C}}

\renewcommand{\rho}{\varrho}

\newcommand{\hil}{\mathcal{H}}

\begin{document}
   
\title[Boundary points of the numerical range]{An observation concerning boundary points of the numerical range}

\author[M. Hansmann]{Marcel Hansmann}
\address{Faculty of Mathematics\\ 
Chemnitz University of Technology\\
Chemnitz\\
Germany.}
\email{marcel.hansmann@mathematik.tu-chemnitz.de}

 \begin{abstract}
A theorem of H\"ubner states that non-round boundary points of the numerical range of a linear operator, i.e. points where the boundary has infinite curvature, are contained in the spectrum of the operator. In this note, answering a question of Salinas and Velasco, we will show that H\"ubner's result remains true under the weaker assumption that the boundary has infinite upper curvature. Our short and simple proof relies on some classical ideas of Berberian.
 \end{abstract}

 
\maketitle 

The numerical range of a bounded linear operator $A$ on a complex Hilbert space $(\hil, \langle . , . \rangle)$ is defined as the set 
$$W(A)= \{ \langle Af,f \rangle : f \in \hil, \|f\| = 1\}.$$
It is well known that the numerical range is a bounded convex set, that it need neither be open nor closed and that its closure contains the spectrum $\sigma(A)$ of $A$. 

An interesting problem is to find geometric conditions on the topological boundary of $W(A)$, which assure that a given point on the boundary is contained in the spectrum of $A$. The most classical result in this direction is due to Donoghue \cite{MR0096127}, who showed that a corner point of $\partial W(A)$, if contained in $W(A)$, is an eigenvalue of $A$. Corner points that are not contained in $W(A)$ are elements of the approximate point spectrum $\sigma_{ap}(A)$, as has been shown by Hildebrandt \cite{MR0200725}. Let us recall that the approximate point spectrum consists of those $\lambda \in \sigma(A)$ for which there exists a sequence of unit vectors $(u_n) \subset \hil$ with $(A-\lambda)u_n \to 0$ for $n \to \infty$. In particular, the set of eigenvalues $\sigma_p(A)$ is contained in $\sigma_{ap}(A)$.

Some decades after the works of Donoghue and Hildebrandt, H\"ubner \cite{MR1371343} generalized their results  considerably by showing that all points where $\partial W(A)$ has infinite curvature (or non-round points, in the terminology of H\"ubner) are contained in $\sigma_{ap}(A)$. More recently, Salinas and Velasco \cite{MR1800238} showed that it is even sufficient to assume that the boundary has unilateral infinite curvature, meaning that either its right- or left-hand curvature at the given point is infinite. An interesting question raised in \cite[Remark 2.6]{MR1800238} was whether one can relax the assumptions even further and to only assume that the \emph{upper} curvature of the boundary is infinite. The main purpose of this note is to answer this question in the affirmative.

Before doing that, let us shortly recall the difference between points of infinite curvature and points of infinite upper curvature (for a more detailed discussion we refer to \cite{Hansmann14}): Let us consider a convex set $\Omega \subset \mr^2$ and let us pick one of its boundary points, which, for simplicity, we assume not to be a corner point of $\Omega$ (meaning that $\Omega$ is not contained in a sector with vertex  at that point and angle less than $\pi$). Using a suitable translation and rotation of $\Omega$, it is no loss to assume that the boundary point we picked is $0$, that $\Omega$ is contained in the closed upper half-plane and that the real line is the unique supporting line for $\Omega$ passing through $0$. Given all these assumptions, we call $0 \in \partial \Omega$ a point of infinite curvature, if
\begin{equation}
  \label{eq:1}
  \gamma_l(0):=\liminf_{\underset{ x \neq 0, (x,y) \in \partial \Omega}{(x,y) \to 0}} \frac{y}{x^2} = \infty.
\end{equation}
By definition, a corner point is a point of infinite curvature. The term $\gamma_l(0)$ is also called the lower curvature of $\partial \Omega$ at $0$.  The upper curvature $\gamma_u(0)$ is defined similarly, but with a limes superior instead of a limes inferior. We call $0$ a point of infinite upper curvature if $\gamma_u(0)=\infty$. Clearly, every point of infinite curvature is a point of infinite upper curvature. However, the simple example where $\Omega$ is the epigraph of the convex function
$$ f : [-1,1] \to \mr_+, \quad f(x) = \left\{
  \begin{array}{cl}
    x^4, & x \leq 0 \\
    x^{3/2}, & x> 0
  \end{array}\right.$$
shows that the converse must not be true. Here $\gamma_u(0)=\infty$, but $\gamma_l(0)=0$. Quite remarkably, the behavior encountered in this example is generic: A result of Zamfirescu \cite{MR593634} states that for most convex bodies and for most of their boundary points $\lambda$ (in each case meaning all except those in a set of first Baire category) one has $\gamma_u(\lambda)=\infty$ and $\gamma_l(\lambda)=0$. 

Let us return to the question of Salinas and Velasco. A first step towards an answer is the following recent result from \cite{Hansmann14} (we remark that in \cite{Hansmann14} H\"ubner's result was extended to unbounded operators).
\begin{thm}\label{thm1}
  Let $\lambda \in \partial W(A) \cap W(A)$ be a point of infinite upper curvature. Then $\lambda$ is an eigenvalue of $A$.
\end{thm}
The proof of this theorem is almost literally the same as Donoghue's proof of the same result for corner points. Indeed, Donoghue argued that if $\lambda=\langle Af,f \rangle$ (where $\|f\|=1$) is not an eigenvalue, then the numerical range of the compression of $A$ to the linear hull of $\{f,Af\}$ is an ellipse, which is contained in $W(A)$, and with $\lambda$ on its boundary. This (and here the assumption that $\lambda$ is a corner point enters) is only possible if the ellipse is degenerated to a line segment with $\lambda$ being one of its end points. But then $\lambda$ must be an eigenvalue of the compression and hence of $A$, a contradiction. Note that Donoghue's proof goes through not only for corner points, but for all points $\lambda \in \partial W(A)$ with the property that there does not exist a non-degenerate ellipse $E$ with $\lambda \in \partial E$ and $\overline{E} \subset \overline{W}(A)$. It is not difficult to see that this property characterizes points of infinite upper curvature (see \cite[Lemma 1]{Hansmann14}).    
 
With the following theorem, we obtain the full answer to the question of Salinas and Velasco.
\begin{thm}\label{thm2}
Let $\lambda \in \partial W(A)$ be a point of infinite upper curvature. Then $\lambda \in \sigma_{ap}(A)$. 
\end{thm}

Given Theorem \ref{thm1}, the proof of Theorem \ref{thm2} is simple. First, we introduce an auxiliary Hilbert space $\hil_0$ and a bounded operator $A_0$ on $\hil_0$, using a classical construction of Berberian \cite{MR0133690}: To this end, let $\operatorname{LIM}$ denote some fixed Banach limit on the space $l^\infty(\mn)$ of complex valued bounded sequences. Let $X$ denote the space of all bounded sequences $x=(x_n)_{n \in \mn} \subset \hil$, equipped with the (possibly degenerate) inner product
$$ \left[ x,y \right] := \underset{n \to \infty}{\operatorname{LIM}} \langle x_n, y_n \rangle.$$
Denoting by $X_0 \subset X$ the subspace of all $x \in X$ with $[x,x]=0$, we define the Hilbert space $\hil_0$ as the completion of the Pre-Hilbert space $(X/X_0, [. , . ])$. Furthermore, we define a bounded operator $A_0$ on $X/X_0$ by setting
$$ A_0((x_n)_{n \in \mn}+X_0) := (Ax_n)_{n \in \mn}+X_0, \qquad (x_n)_{n \in \mn} \in X.$$
One can check that this operator is well defined and continuous with $\|A_0\|=\|A\|$. In particular, there exists a unique continuous extension of $A_0$ to $\hil_0$, which we continue to denote by $A_0$. The purpose of all these constructions are the following two identities, proved by Berberian \cite{MR0133690} and Berberian and Orland \cite{MR0212588}, respectively:
$$\sigma_{ap}(A)=\sigma_p(A_0) \quad \text{and} \quad  \overline{W}(A)=W(A_0).$$ 
In particular, the numerical range of $A_0$ is closed. Now our assumption on $\lambda$ implies that $\lambda \in \partial W(A_0) \cap W(A_0)$ is a point of infinite upper curvature of $\partial W(A_0)$, so by Theorem \ref{thm1} we obtain that $\lambda \in \sigma_p(A_0)=\sigma_{ap}(A)$, which concludes the proof of Theorem \ref{thm2}.    

We conclude this note with a result about points of infinite upper curvature, which are not corner points. Recall that the essential spectrum of $A$ consists of all $\lambda \in \mc$ with the property that $\lambda -A$ is not a Fredholm operator.

\begin{thm}\label{thm3}
Let $\lambda \in \partial W(A)$ be a point of infinite upper curvature, which is not a corner point. Then $\lambda$ is an element of the essential spectrum of $A$. 
\end{thm}      
For points of infinite curvature, this result has been proved independently by Farid \cite{MR1731863}, Salinas and Velasco \cite{MR1800238} and Spitkovsky \cite{MR1752163}. Spitkovsky's proof reduces the problem to H\"ubner's original result about boundary points of infinite curvature. Using Theorem \ref{thm2} instead of H\"ubner's original result in Spitkovsky's argument, we obtain Theorem \ref{thm3}.

 \def\cydot{\leavevmode\raise.4ex\hbox{.}} \def\cprime{$'$}

\end{document}